\documentstyle[11pt]{article}

\begin{document}

\begin{titlepage}

\title{Equivariant Higher Analytic Torsion and Equivariant Euler
Characteristic}

\author{Ulrich Bunke\thanks{Mathematisches Institut, Universit\"at G\"ottingen,
Bunsenstr. 3-5, 37073 G\"ottingen, Germany, bunke@uni-math.gwdg.de}}
\end{titlepage}

\newcommand{\diag}{{\rm diag}}
\newcommand{\proof}{{\it Proof.$\:\:\:\:$}}
\newcommand{\Bbb}{\rm}
\newcommand{\dist}{{\rm dist}}
\newcommand{\kaaa}{{\bf k}}
\newcommand{\paaa}{{\bf p}}
\newcommand{\taaa}{{\bf t}}
\newcommand{\haaa}{{\bf h}}
\newcommand{\R}{{\bf R}}
\newcommand{\Q}{{\bf Q}}
\newcommand{\Z}{{\bf Z}}
\newcommand{\C}{{\bf C}}
\newcommand{\K}{{\tt K}}
\newcommand{\Naaa}{{\bf N}}
\newcommand{\gaaa}{{\bf g}}
\newcommand{\maaa}{{\bf m}}
\newcommand{\aaaa}{{\bf a}}
\newcommand{\naaa}{{\bf n}}
\newcommand{\brr}{{\bf r}}
\newcommand{\res}{{\rm res}}
\newcommand{\Tr}{{\rm Tr}}
\newcommand{\cT}{{\cal T}}
\newcommand{\dom}{{\rm dom}}
\newcommand{\db}{{\bar{\partial}}}
\newcommand{\g}{{\gaaa}}
\newcommand{\cZ}{{\cal Z}}
\newcommand{\cH}{{\cal H}}
\newcommand{\cM}{{\cal M}}
\newcommand{\interi}{{\rm int}}
\newcommand{\singsupp}{{\rm singsupp}}
\newcommand{\cE}{{\cal E}}
\newcommand{\cV}{{\cal V}}
\newcommand{\cI}{{\cal I}}
\newcommand{\cC}{{\cal C}}
\newcommand{\mod}{{\rm mod}}
\newcommand{\cK}{{\cal K}}
\newcommand{\cA}{{\cal A}}
\newcommand{\cEp}{{{\cal E}^\prime}}
\newcommand{\cU}{{\cal U}}
\newcommand{\Hom}{{\mbox{\rm Hom}}}
\newcommand{\vol}{{\rm vol}}
\newcommand{\cO}{{\cal O}}
\newcommand{\End}{{\mbox{\rm End}}}
\newcommand{\Ext}{{\mbox{\rm Ext}}}
\newcommand{\rk}{{\mbox{\rm rank}}}
\newcommand{\im}{{\mbox{\rm im}}}
\newcommand{\sign}{{\rm sign}}
\newcommand{\spann}{{\mbox{\rm span}}}
\newcommand{\symm}{{\mbox{\rm symm}}}
\newcommand{\cF}{{\cal F}}
\newcommand{\Ree}{{\rm Re }}
\newcommand{\Res}{{\mbox{\rm Res}}}
\newcommand{\Imm}{{\rm Im}}
\newcommand{\inter}{{\rm int}}
\newcommand{\clo}{{\rm clo}}
\newcommand{\tg}{{\rm tg}}
\newcommand{\ee}{{\rm e}}
\newcommand{\Li}{{\rm Li}}
\newcommand{\cN}{{\cal N}}
 \newcommand{\conv}{{\rm conv}}
\newcommand{\op}{{\mbox{\rm Op}}}
\newcommand{\tr}{{\mbox{\rm tr}}}
\newcommand{\cs}{{c_\sigma}}
\newcommand{\ctg}{{\rm ctg}}
\newcommand{\degg}{{\mbox{\rm deg}}}
\newcommand{\Ad}{{\mbox{\rm Ad}}}
\newcommand{\ad}{{\mbox{\rm ad}}}
\newcommand{\codim}{{\mbox{\rm codim}}}
\newcommand{\Gr}{{\mbox{\rm Gr}}}
\newcommand{\coker}{{\rm coker}}
\newcommand{\id}{{\mbox{\rm id}}}
\newcommand{\ord}{{\mbox{\rm ord}}}
\newcommand{\nat}{{\bf  N}}
\newcommand{\supp}{{\mbox{\rm supp}}}
\newcommand{\spec}{{\mbox{\rm spec}}}
\newcommand{\Ann}{{\mbox{\rm Ann}}}
\newcommand{\aca}{{\aaaa_\C^\ast}}
\newcommand{\acag}{{\aaaa_{\C,good}^\ast}}
\newcommand{\acage}{{\aaaa_{\C,good}^{\ast,extended}}}
\newcommand{\ck}{{\cal K}}
\newcommand{\tck}{{\tilde{\ck}}}
\newcommand{\tnk}{{\tilde{\ck}_0}}
\newcommand{\ceep}{{{\cal E}(E)^\prime}}
 \newcommand{\ncE}{{{}^\naaa\cE}}
 \newcommand{\Or}{{\rm Or }}

\newcommand{\cB}{{\cal B}}
\newcommand{\hc}{{{\cal HC}(\gaaa,K)}}
\newcommand{\hcma}{{{\cal HC}(\maaa_P\oplus\aaaa_P,K_P)}}

\newcommand{\vsl}{{V_{\sigma_\lambda}}}
\newcommand{\czg}{{\cZ(\gaaa)}}
\newcommand{\csl}{{\chi_{\sigma,\lambda}}}
\newcommand{\cR}{{\cal R}}
\def\hB{\hspace*{\fill}$\Box$ \newline\noindent}
\newcommand{\varho}{\varrho}
\newcommand{\ind}{{\rm index}}
\newcommand{\Ind}{{\rm Ind}}
\newcommand{\Fin}{{\mbox{\rm Fin}}}
\newcommand{\cS}{{\cal S}}
\newcommand{\orig}{{\cal O}}
\def\hB{\hspace*{\fill}$\Box$ \\[0.5cm]\noindent}
\newcommand{\cL}{{\cal L}}
 \newcommand{\cG}{{\cal G}}
\newcommand{\npci}{{\naaa_P\hspace{-1.5mm}-\hspace{-2mm}\mbox{\rm coinv}}}
\newcommand{\pki}{{(\paaa,K_P)\hspace{-1.5mm}-\hspace{-2mm}\mbox{\rm inv}}}
\newcommand{\mki}{{(\maaa_P\oplus \aaaa_P, K_P)\hspace{-1.5mm}-\hspace{-2mm}\mbox{\rm inv}}}

\newcommand{\npi}{{\naaa_P\hspace{-1.5mm}-\hspace{-2mm}\mbox{\rm inv}}}
\newcommand{\ngp}{{N_\Gamma(\pi)}}
\newcommand{\gbg}{{\Gamma\backslash G}}
\newcommand{\gkm}{{ Mod(\gaaa,K) }}
\newcommand{\ggkm}{{  (\gaaa,K) }}
\newcommand{\pkm}{{ Mod(\paaa,K_P)}}
\newcommand{\ppkm}{{  (\paaa,K_P)}}
\newcommand{\makm}{{Mod(\maaa_P\oplus\aaaa_P,K_P)}}
\newcommand{\mmakm}{{ (\maaa_P\oplus\aaaa_P,K_P)}}
\newcommand{\cP}{{\cal P}}
\newcommand{\gm}{{Mod(G)}}
\newcommand{\gk}{{\Gamma_K}}
\newcommand{\La}{{\cal L}}
\newcommand{\cug}{{\cU(\gaaa)}}
\newcommand{\cuk}{{\cU(\kaaa)}}
\newcommand{\dc}{{C^{-\infty}_c(G) }}
\newcommand{\gdk}{{\gaaa/\kaaa}}
\newcommand{\dgkm}{{ D^+(\gaaa,K)-\mbox{\rm mod}}}
\newcommand{\dgm}{{D^+G-\mbox{\rm mod}}}
\newcommand{\vect}{{\C-\mbox{\rm vect}}}
 \newcommand{\cig}{{C^{ \infty}(G)_{K} }}
\newcommand{\gami}{{\Gamma\hspace{-1.5mm}-\hspace{-2mm}\mbox{\rm inv}}}
\newcommand{\cQ}{{\cal Q}}
\newcommand{\mmap}{{Mod(M_PA_P)}}

\newtheorem{prop}{Proposition}[section]
\newtheorem{lem}[prop]{Lemma}
\newtheorem{ddd}[prop]{Definition}
\newtheorem{theorem}[prop]{Theorem}
\newtheorem{kor}[prop]{Corollary}
\newtheorem{ass}[prop]{Assumption}
\newtheorem{con}[prop]{Conjecture}
\newtheorem{prob}[prop]{Problem}
\newtheorem{fact}[prop]{Fact}

\maketitle

\begin{abstract}
We show that J. Lott's equivariant higher analytic torsion for compact group actions
depends only on the equivariant Euler characteristic.
 \end{abstract}

\tableofcontents

\section{Introduction}

Let $G$ be a compact connected Lie group with Lie algebra $g$.
Let $I(G)$ denote the ring of $\Ad(G)$-invariant polynomials
on $g$. Then $I(G)^1:=\{f\in I(G)\:|\:f(0)=0\}$ is a maximal
ideal  of $I(G)$. By $\hat{I}(G)$ we denote the $I(G)^1$-adic
completion of $I(G)$. We define $\tilde{I}(G):=\hat{I}(G)/\C 1$.

Let $M$ be a closed oriented $G$-manifold. Then Lott \cite{lott94}
defined equivariant higher analytic torsion $T(M)$ of $M$ (see Def. \ref{tvonm}).
To be precise, in \cite{lott94}, Def. 2, he defined an element
$T(M,g^M,F)\in \hat{I}(g)$, where $g^M$ is a $G$-equivariant
Riemannian metric and $F$ is a equivariant flat hermitean
vector bundle with trivial momentum map \cite{lott94} (14).
In our case for $F$ we take the trivial flat hermitean bundle
$F:=M\times \C$, where $G$ acts on
the first factor. By \cite{lott94},  Cor. 1, the class
$T(M):=[T(M,g^M,M\times \C)]\in \tilde{I}(g)$ is independent of $g^M$.
By definition $T(M)$ is a differential topological invariant of the
$G$-manifold $M$. If $M$ is even-dimensional, then by \cite{lott94}, Prop. 9,
we have $T(M)=0$.

Let $\Or(G)$ denote the orbit category of  $G$ (see L\"uck \cite{lueck89},
Def. 8.16), and let $U(G)$ be the Euler ring of $G$ (\cite{lueck89}, Def. 5.10).
By \cite{lueck89}, Prop. 5.13, we can identify
$$U(G)=\prod_{[G/H]\in \Or(G)} \Z [G/H]\ ,$$
where the product runs over all isomorphism classes of objects of $\Or(G)$.
If $X$ is a $G$-space of the $G$-homotopy type of a finite $G$-CW complex,
then we can define its equivariant Euler characteristic $\chi_G(X)\in U(G)$.
If ${E_\alpha}$ is the finite collection of $G$-cells of $X$, then
$$\chi_G(X):=\sum_\alpha (-1)^{\dim(E_\alpha)} [G/t(E_\alpha)]\ ,$$
where $t(E)=H$ is the type of the cell $E=G/H\times D^{\dim(E_\alpha)}$
(see \cite{lueck89}, Lemma 5.6).
Any compact $G$-manifold has the $G$-homotopy type of a finite
$G$-CW complex (\cite{lueck89}, 4.36), and thus $\chi_G(M)\in U(G)$
is well defined.

In the present note we define a homomorphism $T_G:U(G)\rightarrow \tilde{I}(G)$
(Lemma \ref{tge}), such that our main result can be formulated as follows.
\begin{theorem}\label{main}
Let $G$ be a compact connected Lie group. If $M$ is a closed oriented
$G$-manifold, then
$$T(M)=T_G \: \chi_G(M)\ .$$
\end{theorem}
This theorem answers essentially the question posed by Lott \cite{lott94}, Note 4.
As we will see below it can be employed to compute $T(M)$ effectively.

Let $H\subset G$ be a closed subgroup. Then by \cite{lueck89}, 7.25 and 7.27,
there is a restriction map $\res^G_H:U(G)\rightarrow U(H)$ such that
$\res^G_H \chi_G(M)=\chi_H(\res^G_H\:M)$ for any compact $G$-manifold,
where $\res^G_H(M)$ denotes $M$ with the induced action of $H$.

The inclusion $h\hookrightarrow g$ induces a map
$\res^G_H:\tilde{I}(G)\rightarrow \tilde{I}(H)$. It is an immediate consequence
of the Definition \ref{tvonm} of $T(M)$, that
\begin{equation}\label{tz1}\res^G_H T(M):= T(\res^G_H \: M)\ .\end{equation}
This is compatible with
\begin{equation}\label{tz2}
\res^G_H \circ T_G = T_H\circ \res^G_H \ .
\end{equation}

Let $S(G)\subset \Or(G)$ be the full subcategory with objects $G/H$, where
$H$ is isomorphic to $S^1$. By Corollary \ref{ress1} the collection
$\res^G_H T(M)$, $G/H\in S(G)$, determines $T(M)$.
In order to compute $T(M)$ it is thus sufficient to define
$T_{S^1}:U(S^1)\rightarrow \tilde{I}(S^1)$. If $H\subset G$ is
isomorphic to $S^1$, then $T_H$ is defined, and we have
$$\res^G_H T(M)=T(\res^G_H M)=T_H \chi_H \res^G_H(M)=T_H \res^G_H \chi_G(M)\ .$$
In order to give  an  explicit formula for $T(M)$ in terms of the $G$-homotopy
type of $M$ it remains to give the formula for $T_{S^1}$.

Since $T_{S^1}$ has to satisfy Theorem \ref{main}, we are forced to put
\begin{eqnarray}\label{ts1}
T_{S^1}([S^1/S^1])&=&T(*)=0\\
T_{S^1}([S^1/H])&=&T(S^1/H), \quad H\not=S^1\nonumber\ .
\end{eqnarray}
For $n\in \nat$ let $F_n:S^1\rightarrow S^1$ be 
the $n$-fold covering. The derivative $F_{n*}$ of $F_n$ at $1\in S^1$ is
multiplication by $n$.
By $\tilde{F}_n:\tilde{I}(S^1)\rightarrow \tilde{I}(S^1)$ we denote 
the induced map.
If $H\subset S^1$ is different from $S^1$, then it is a cyclic subgroup
of finite order $|H|$. It is again an easy consequence of the Definition \ref{tvonm}
of $T(M)$, that 
\begin{equation} \label{uy1} T(S^1/H)=\tilde{F}_{|H|} T(S^1)\ .\end{equation}

Let $S^1:=\{z\in \C\:|\: |z|=1\}$. We identify $s^1\cong \R$ such that
the  exponential map is given by $\exp(y):=\ee^{iy}$. Then $I(S^1)=\C[y]$, and
we identify $\tilde{I}(S^1)\cong y\C(y)$.
By \cite{lott94}, Prop.  11,  we then have
$$T(S^1)=2\sum_{k=1}^\infty \left(\begin{array}{c}4k\\2k\end{array}\right) \Li_{2k+1}(1)\left(\frac{y}{8\pi}\right)^{2k}\ ,$$
where $$\Li_j(z):=\sum_{m=1}^\infty \frac{z^m}{m^j}\ .$$
It follows that
$$T(S^1/H)=2\sum_{k=1}^\infty \left(\begin{array}{c}4k\\2k\end{array}\right) \Li_{2k+1}(1)\left(\frac{y|H|}{8\pi}\right)^{2k}\ .$$

We now discuss some consequences.

\begin{lem}
If $M$ is a closed oriented $S^1$-manifold, then
$T(M)$  and $\chi(M)$  together determine $\chi_{S^1}(M)$.
\end{lem}
\proof
$\chi(M)$ is the coefficient at $[S^1/S^1]$ of $\chi_{S^1}(M)$.
Let $\{H_1, \dots, H_l\}$ be the finite set of orbit types of $M$ with
$H_i\not=S^1$. Since $\Li_{j}(1)\not=0$ for all $j\in \nat$, $j\ge 2$,
the torsion $T(M)$ determines the numbers $r_l:=\sum_{i=1}^l |H_i|^j$, $j\in 2\nat$.
But  vice versa the numbers $r_l$ determine $|H_i|$ and therefore $H_i$, $i=1,\dots l$.\hB

\begin{lem}\label{torus}
Let $T$ be a $k$-dimensional torus and $H\subset T$ be a closed 
subgroup. If $\dim(T/H)\ge 2$, then $T(T/H)=0$, and if $\dim(T/H)=1$,
then $T(T/H)=\tilde{P}(T(S^1))$, where
$\tilde{P}:\tilde{I}(S^1)\rightarrow \tilde{I}(T)$ is induced
by the projection $P:T\rightarrow T/H\cong S^1$.
\end{lem}
\proof
Let $R\in S(T)$. Then $\chi_R (\res^T_R T/H)=\chi((T/H)/R)[R/R\cap H]$.
If $\dim(T/H)\ge 2$, then $(T/H)/R$ is a torus and $\chi((T/H)/R)=0$.
If $\dim(T/H)=1$, then $\chi((T/H)/R)\not=0$ iff $(T/H)/R$ is a point.
Thus $\chi_R(\res^T_R T/H)=[R/R\cap H]$.
The Lemma now follows from  (\ref{tz2}) and (\ref{uy1}).\hB

Let $G/K$ be a compact symmetric space associated to the Cartan involution
$\theta$ of $G$. We fix a $\theta$-stable maximal torus $T\subset G$.
Then $T\cap K=:S$ is a maximal compact torus of $K$. The rank of
$G/K$ is by definition $\rk(G/K):=\dim(T)-\dim(S)$.
Let $W_G(T)$, and $W_K(T)$ be the Weyl groups
of $(G,T)$ and $(K,T)$. If $\rk(G/K)=1$, then for $w\in W_G(T)$ we have a projection
$P_w:T\rightarrow T/S^w\cong S^1$, where $S^w=wSw^{-1}$. It induces
a map $\tilde{P}_w:\tilde{I}(S^1)\rightarrow \tilde{I}(T)$.
Since $\res^G_T:\tilde{I}(G)\rightarrow \tilde{I}(T)$ is injective,
the following Lemma gives an explicit computation of $T(G/K)$.
\begin{lem}
If $\rk(G/K)\ge 2$, then $T(G/K)=0$, and if $\rk(G/K)=1$, then
$\res^G_T T(M)=\sum_{W_G(T)/W_K(T)} \tilde{P}_w(T(S^1))$.
\end{lem}
\proof
Fix $S^1\cong R\subset T$.
If $H\subset T$ is a closed subgroup, then $\chi_R(\res^T_R T/H)=0$ except
if $\dim(T/H)=1$.
In \cite{bunke972} we have shown that
$$\chi_T(\res^G_T \:G/K)=  \sum_{W_G(T)/W_K(T)} [T/S^w] + \mbox{higher dimensional staff}\ .$$
Hence if $\rk(G/K)\ge 2$, then by Lemma \ref{torus} $\res^G_T T(M)=0$, and if $\rk(G/K)=1$, then
$$\res^G_T T(G/K)=  \sum_{W_G(T)/W_K(T)} T[T/S^w]\ .$$
Applying \ref{torus} we obtain the desired result.
\hB

We now briefly describe the contents of the remainder of the paper.
In Section \ref{sumform} we prove our main analytic result
Theorem \ref{sumff} saying that $T(M)$ is essentially additive.
In Section \ref{covprod} we study the behaviour of $T(M)$ under coverings
and with respect to cartesian products.
In Section \ref{formalities} we extend the analytic results to manifolds
with corner singularities using certain formal considerations.
In Section \ref{subgr} we show
that $T(M)$ is determined by its restrictions to all subgroups $H\cong S^1$.
In Section \ref{s1}  we first prove Theorem \ref{main} for $G=S^1$, and then
we construct $T_G$ and finish the proof of Theorem \ref{main} for general $G$.

\section{Additivity of equivariant torsion}\label{sumform}

We first recall the definition of higher equivariant torsion
\cite{lott94}, Def. 2. Let $G$ be a connected Lie group
with Lie algebra $g$. Let $M$ be a closed oriented $G$-manifold.
We write $\Omega(M):=C^\infty(M,\Lambda^*T^*M)$
and  $d:\Omega(M)\rightarrow \Omega(M)$ for the
differential of the de Rham complex.

For $X\in g$ let $X^*\in C^\infty(M,TM)$ denote the corresponding
fundamental vector field. We set
$$I:=\sum_{\alpha} X^\alpha\otimes i_{X^*_\alpha}\in S(g^*)\otimes \End(\Omega(M))\ ,$$
where $X_\alpha\in g$, $X^\alpha\in g^*$ run over a base of $g$ or dual base
of $g^*$, respectively, and $i_Y$ denotes interior
multiplication by the vector field $Y$.
We choose a $G$-invariant Riemannian metric $g^M$.
It induces a pre Hilbert space structure on $\Omega(M)$,
and we let $e_Y$ be the adjoint of $i_Y$.
We set $E:=\sum_\alpha X^\alpha\otimes e_{X^*_\alpha}$.

For $t>0$ we define
$$d_t:=\sqrt{t}d-\frac{1}{4\sqrt{t}}I,\quad
\delta_t:= \sqrt{t} d^*+\frac{1}{4\sqrt{t}}E\ .$$
Then we put
\begin{equation}\label{hamil}D_t:=\delta_t-d_t\in S(g^*)\otimes \End(\Omega(M))\ .\end{equation}
Let $S(g^*)^1:=\{f\in S(g^*)\:|\: f(0)=0\}$, and let $\hat{S}(g^*)$
be the $S(g^*)^1$-adic completion. Since
$$D_t^2=-t  \Delta \quad (\mbox{mod $S^1(g^*)\otimes \End(\Omega(M))$})$$
we can form
$$\ee^{D_t^2} \in \hat{S}(g^*)\otimes \End(\Omega(M))\ .$$
Moreover we have
$$\Tr_s N \ee^{D_t^2}\in \hat{I}(G)\ ,$$
where $N$ is the $\Z$-grading operator on $\Omega(M)$, and
$\Tr_s$ is the  $\Z_2$-graded trace on $\End(\Omega(M))$.
Define $\chi^\prime(M):=\sum_{p=0}^\infty p (-1)^p \dim\: H^*(M,\R)$.
Then the function
$$s\mapsto -\frac{1}{\Gamma(s)}\int_0^\infty (\Tr_s N \ee^{D_t^2}-\chi^\prime(M)) t^{s-1} dt$$
is holomorphic for $\Ree(s)>>0$, and it has a meromorphic continuation
to all of $\C$ which is regular at $s=0$.

\begin{ddd}\label{tvonm}
The equivariant higher torsion $T(M)\in\tilde{I}(G)$ of the
$G$-manifold $M$ is represented by
$$-\frac{d}{ds}_{|s=0}\frac{1}{\Gamma(s)}\int_0^\infty (\Tr_s N \ee^{D_t^2} -\chi^\prime(M)) t^{s-1} dt\ .$$
\end{ddd}

If $M$ is odd-dimensional, then by \cite{lott94}, Cor. 1, $T(M)$ is independent of the choice of the $G$-invariant
Riemannian metric $g^M$.
If $M$ is even-dimensional, the by \cite{lott94}, Prop 9,  
we have $T(M)=0$.

Let $M$ be a closed oriented $G$-manifold, and let $N$ be a
$G$-invariant oriented hypersurface such that $M\setminus N$ has two components, i.e. there
are compact manifolds $M_1$, $M_2$ with boundary $\partial M_i=N$, $i=1,2$
such that $M=M_1\cup_N M_2$.
We form the closed oriented $G$-manifolds $\tilde{M_i}:=M_i\cup_NM_i$, the doubles
of $M_i$.

\begin{theorem}\label{sumff}
$$2T(M)=T(\tilde{M}_1)+T(\tilde{M}_2)\ .$$
\end{theorem}
\proof
We choose Riemannian metrics on $M$ and $\tilde{M}_i$, $i=1,2$.
Then let $D_t$ and $D_{t,i}$, $i=1,2$ denote the operators (\ref{hamil}) for
$M$ and $\tilde{M}_i$, respectively.
We define $\delta(t)\in \hat{I}(G)$ by
$$\delta(t):=2\Tr_s N\ee^{D_t^2}-\Tr_s  N\ee^{D_{t,1}^2}- \Tr_s  N\ee^{D_{t,2}^2}
 -( 2\chi^\prime(M)-\chi^\prime(\tilde{M}_1)-\chi^\prime(\tilde{M}_2))\ .$$
We have to show that
$$0=[-\frac{d}{ds}_{|s=0}\frac{1}{\Gamma(s)} \int_0^\infty \delta(t) t^{s-1}   dt]  \ ,$$
where $[.]$ denotes the class of $"."$ in $\tilde{I}(G)$.

We now specialize the choice of Riemannian metrics.
We choose a $G$-invariant collar neighbourhood $(-1,1)\times N\hookrightarrow M$
such that $\{0\}\times N$ is mapped to $N$.
Then we assume that $g^M$ is a product metric $dr^2+g^N$ on the collar.
The metric $g^M$ induces natural Riemannian metrics $g^{\tilde{M}_i}$ on $\tilde{M}_i$.

For $R>1$ let $g^M(R)$ be the Riemannian metric which coincides
with $g^M$ outside the collar, and which is such that the collar
is isometric to $(-R,R)\times N$.
Similarly we obtain metrics $g^{\tilde{M}_i}(R)$ on $\tilde{M}_i$.

Let $\delta(t,R)$ be defined with respect to these choices of metrics.
While $\delta(t,R)$ may depend on $R$, it is known that
$$[-\frac{d}{ds}_{|s=0}\frac{1}{\Gamma(s)} \int_0^\infty \delta(t,R) t^{s-1}   dt]\in \tilde{I}(G)$$
is independent of $R$. The proof of the theorem is obtained
by studying the behaviour of $\delta(t,R)$ as $R$ tends to infinity.

Note that $\hat{I}(G)$ is a locally convex topological vector space.
\begin{prop}\label{pr1}
For any seminorm $|.|$ on $\hat{I}(G)$ 
there are constants $C<\infty,c>0$ such that for all $t>0$, $R>1$
$$|\delta(t,R)|<C\ee^{-\frac{c R^2}{t}}\ .$$
\end{prop}
\proof
This follows from a standard argument using the finite propagation speed
method \cite{cheegergromovtaylor82}. We leave the details to the interested reader.
\hB

Let $I(G)_1\subset \hat{I}(G)$ be the closed subspace of at most
linear invariant polynomials on $g$ and put $\check{I}(G):=\hat{I}(G)/I(G)_1$.
By $[[.]]$ we denote classes in this topological quotient space.

\begin{prop}\label{pr2}
For any seminorm  $|.|$ on $\check{I}(G)$
there is a constant $C<\infty$ such that for all $R>1$, $t>1$
$$|[[\delta(t,R)]]| < C t^{-1} R \  .$$
\end{prop}
\proof
This is a consequence of the more general estimate
\begin{equation}\label{est}|[[\Tr_s N\ee^{D_t(R)^2}]]|< C t^{-1} R \end{equation}
which also holds for $M$ replaced by $\tilde{M}_i$.
Here $D_t(R)$ denotes the operator (\ref{hamil}) associated to $g^M(R)$.

We can assume that $|.|$ is the restriction to $\check{I}(G)$ of a seminorm of $\hat{S}(g^*)/S_1(g^*)$, where $S_1(g^*)$ denotes the subspace $\C\oplus g^*$.
There is an $m>0$ depending on $|.|$ such that $|[[U]]|=0$ for all $U\in \hat{S}(\gaaa^*)^{m}$.
Let $\Delta(R)$ denote the Laplace operator on
differential forms associated to the Riemannian metric $g^M(R)$.
We have $$D_t^2(R)=-t\Delta(R) +  \cN + \frac{1}{t} \cN_1\ ,$$
(to be precise we should write $\cN(R),\cN_1(R)$) where
\begin{eqnarray*}
\cN&:=&\frac{1}{4}[d^*(R)-d,E+I]\\
\cN_1&:=& Q \\
Q&:=& \frac{1}{16} [I,E]\ ,
\end{eqnarray*}
(the commutators are understood in the graded sense)
belong to $S(\gaaa^*)^1\otimes \End(\Omega(M))$.
  
As in \cite{berlinegetzlervergne92}, 9.46,  we write
\begin{eqnarray}
\Tr_s N\ee^{D_t(R)^2}&=&\sum_{k=0}^\infty  \int_{\Delta_k} U_k(\sigma,R) d\sigma   \ , \label{volt}             \\
U_k(\sigma,R)&:=&\Tr_s N \ee^{-\frac{1}{4}t\sigma_0 D(R)^2}(\cN + \frac{1}{t} \cN_1)\dots (\cN + \frac{1}{t} \cN_1)\ee^{-\frac{1}{4}t\sigma_k D(R)^2}\ ,\nonumber
\end{eqnarray}
where $\Delta_k\subset \R^{k+1}$ denotes the standard simplex
such that $\Delta_k\ni \sigma=(\sigma_0,\dots,\sigma_k)$ satisfies 
$\sum_{i=0}^k \sigma_i=1$.

The Riemannian metric $\g^M(R)$ induces a pre Hilbert space structure on $\Omega(M)$.
The trace (operator) norm $\|.\|_1$ ($\|.\|$) on $\End(\Omega(M))$ and $|.|$
together induce norms on $\hat{S}(g)/S_1(g)\otimes \End(\Omega(M))$ which we also denote by $\|.\|_1$ ($\|.\|)$).

\begin{lem}\label{trace}
There is a constant $C<\infty$ such that for all $t>1$ and $R>1$ we have
$$\|\ee^{-t\Delta(R)}\|_1 < C R\ .$$
\end{lem}
\proof
The operator $\ee^{-t\Delta(R)}$ is positive. Thus
$\|\ee^{-t\Delta(R)}\|_1=\Tr\: \ee^{-t\Delta(R)}$.
Let $W(t,x,y)(R)$ be the integral kernel of $\ee^{-t\Delta(R)}$.
The family $(M,g^M(R))$ of Riemannian manifolds
has uniformly bounded geometry as $R$ varies in $[1,\infty)$,
i.e. there are uniform curvature bounds, and the injectivity radius
is uniformly bounded from below.
Standard heat kernel estimates (see e.g.
\cite{cheegergromovtaylor82})
imply that there is a constant $C_1<\infty$ such that
for all $x\in M$, $t>1$, $R>1$ we have $|W(t,x,x)(R)| < C_1$.
In particular, for some $C,C_2<\infty$ independent of $R>1$, $t>1$ we have
\begin{eqnarray*}
\Tr\: \ee^{-t\Delta(R)}&=& \int_M \tr\: W(t,x,x)(R) \vol_{g^M(R)}(x)\\
&<& C_2 \vol_{g^M(R)}(M)\\
&<& C R \ .
\end{eqnarray*}
This finishes the proof of the lemma. \hB

\begin{lem}\label{ooo}
There is a $C<\infty$ such that for all $R>1$ and $t,s> 0$ we have
$$\|[[\ee^{-t\Delta(R)}\cN\ee^{-s\Delta(R)}]]\| < C (t^{-1/2}+s^{-1/2})\ .$$
\end{lem}
\proof
Since $\cN=[d^*(R)-d,E+I]$ and
$\|E+I\|$ is uniformly bounded w.r.t. $R$
it suffices to show that
there exists $C_1<\infty$ such that for all $R>1$ and $t>0$ we have
$$\|\ee^{-t\Delta(R)}d\|<C_1t^{-1/2},\quad \|\ee^{-t\Delta(R)}d^*(R)\|<C_1t^{-1/2}\ .$$
We consider the first estimate.
Note that $dd^*(R) + d^*(R)d=\Delta(R)$, and the ranges of
$dd^*(R)$ and $d^*(R)d=$ are perpendicular. Thus
\begin{eqnarray*}
\|\ee^{-t\Delta(R)}d\|&=&\|\ee^{-t\Delta(R)}dd^*(R)\ee^{-t\Delta(R)}\|^{1/2}\\
&\le&\|\ee^{-t\Delta(R)}\Delta(R)\ee^{-t\Delta(R)}\|^{1/2}\\
&=&t^{-1/2} \|\ee^{-t\Delta(R)}t\Delta(R)\ee^{-t\Delta(R)}\|^{1/2}\\
&\le&t^{-1/2} \:\sup_{x\ge 0} x\ee^{-x}\\
&\le&C_1 t^{-1/2}\ .
\end{eqnarray*}
\hB
If $A$ is of trace class and $B$ is bounded, then we have
$|\Tr \: AB|\le \|B\| \|A\|_1$. Note that $\|\cN_1\|$ is uniformly
bounded w.r.t. $R$.
Applying this and Lemmas \ref{ooo} and \ref{trace} to $U_k$
we obtain $C, C_1 <\infty$ such that for all $R>1$ and $t>1$
we  have
\begin{eqnarray}|[[U_k(\sigma,R)]]| &<& C_1 R t^{-k/2} \sum_{i=0}^k \sigma_i^{-1/2} \nonumber\\
|[[\int_{\Delta_k} U_k(\sigma,R) d\sigma]]|  &<& C t^{-k/2} R \ .\label{tz4}\end{eqnarray}
Note that $|[[U_k]]|=0$ for $k>m$. In order to obtain
(\ref{est}) from (\ref{volt}) and (\ref{tz4})
it remains to discuss $U_1$. Since $\cN\in S(g)^1$ there exists
$C,C_1<\infty$ such that for all $R>1$ and $t>1$
\begin{eqnarray*}
|[[U_1(\sigma,R)]]|&=& |[[\Tr_s N \ee^{-t\sigma_0 \Delta(R)}\frac{1}{t} \cN_1 \ee^{-t\sigma_1 \Delta(R)}]]|\\
&=& |[[\Tr_s N  \frac{1}{t} \cN_1 \ee^{-t \Delta(R)}]]| \\
&<& C_1 R t^{-1}\\
|[[\int_{\Delta_1} U_1(\sigma,R) d\sigma]]| &< & C  R t^{-1}
\end{eqnarray*}
This finishes the proof of the proposition.\hB

We now continue with the proof of the theorem.
Let $|.|$ any seminorm on $\check{I}(G)$ as in the proof of Proposition \ref{pr2}.
By Propositions \ref{pr1} and \ref{pr2} we can write
$$\sigma(R):=-\frac{d}{ds}_{|s=0}\frac{1}{\Gamma(s)} \int_0^\infty [[\delta(t,R)]] t^{s-1}   dt\ ,$$
and the integral converges at $t=0$ and $t=\infty$ uniformly in $s\in (-1/2,1/2)$.
We can perform the derivative  and  obtain
\begin{eqnarray*}\sigma(R) &=&-\int_0^\infty [[\delta(t,R)]] t^{-1} dt  \\
&=&-\int_0^R [[\delta(t,R)]] t^{-1} dt + \int_R^\infty [[\delta(t,R)]] t^{-1} dt  \ .
\end{eqnarray*}
By Proposition \ref{pr1} there are $C_1<\infty$, $c_1>0$ such that
for all $R>1$ we have
\begin{eqnarray*} |[[\int_0^{R^{3/2}} \delta(t,R) t^{-1} dt]]| &\le&
                            \int_0^{R^{3/2}} C \ee^{-\frac{c R^2}{t}} t^{-1} dt \\
&\le & C_1 \ee^{-c_1 R^{1/2}}\ .
\end{eqnarray*}
Moreover by Proposition \ref{pr2} there is a $C<\infty$ such that for all $R>1$
\begin{eqnarray*} |[[\int_{R^{3/2}}^\infty \delta(t,R) t^{-1} dt]]| &\le&  \int_{R^{3/2}}^\infty C R t^{-2} dt\\
&=& C  R^{-1/2} \ .
\end{eqnarray*}
We now let $R$ tend to infinity and take into account that
$\sigma(R)$ is independent of $R$ in order to conclude
that \begin{equation}\label{rrr1}\sigma(R)=0\ .\end{equation}
We have shown that $[[T(M)]]=[[T(\tilde{M}_1)]]+[[T(\tilde{M}_2)]]$. 

We now  consider the remaining component $T_1(M)\in I_1(G)/\C 1$.
Note that $\cN=-\frac{1}{2}L + [d,E]+[d^*,I]$, where
$L:=\sum_\alpha  X^\alpha\otimes L_{X_\alpha^*}$ and
$L_Y$ denotes the Lie derivative with respect to the vector field $Y$.
Since $[d,E],$ and $[d^*,I]$ shift the form degree by $\pm 2$ we obtain
$$\Tr_s N \cN \ee^{-t \Delta}= - \frac{1}{2}\Tr_s N L \ee^{-t \Delta}\ .$$

Let $\rho_{an}(M,g^M):G\rightarrow \C$
denote the equivariant analytic torsion defined by \cite{lottrothenberg91}
$$\rho_{an}(M,g^M)(g):=
-\frac{d}{ds}_{|s=0}\frac{1}{\Gamma(s)}\int_0^\infty (\Tr_s N g\ee^{t\Delta}
-\chi^\prime(M)) t^{s-1} dt \ .$$
If we define 
$$\delta(t,g):=  2\Tr_s N g\ee^{t\Delta}-\Tr_s  N g\ee^{t\Delta_1}- \Tr_s  N g\ee^{t\Delta_2}
 -( 2\chi^\prime(M)-\chi^\prime(\tilde{M}_1)-\chi^\prime(\tilde{M}_2))\ ,
$$
then    there are  $C<\infty$, $c>0$ such that for all $g\in G$
\begin{eqnarray*}
|\delta(t,g)|\le C \ee^{-\frac{c}{t}}&&\forall t\in (0,1]\\
|\delta(t,g)|\le C\ee^{-ct}&&\forall t\in [1,\infty)\ .\end{eqnarray*}
The first estimate is again a consequence of the finite propagation speed
method \cite{cheegergromovtaylor82}.
Similar estimates hold for the derivative of $\delta(t,g)$ w.r.t. $g$. 
We have 
$$\sigma_1(g):=-\int_0^\infty \delta(t,g) t^{-1} dt=2\rho_{an}(M,g^M)-\rho_{an}(\tilde{M}_1,g^{\tilde{M}_1})-\rho_{an}(\tilde{M}_2,g^{\tilde{M}_2})\ .$$
On the one hand in \cite{bunke972} we have shown that on the dense
subset of $G$ consisting of elements of finite order 
$$2\rho_{an}(M,g^M)-\rho_{an}(\tilde{M}_1,g^{\tilde{M}_1})-\rho_{an}(\tilde{M}_2,g^{\tilde{M}_2}) = const\ .$$
On the other hand $\sigma_1$ is differentiable.
We conclude
\begin{eqnarray*}
0&=&d_{|g=1} \sigma_1\\
&=&-\int_0^\infty d_{|g=1}\delta(t,.) t^{-1} dt\\
&=&-\int_0^\infty (2\Tr_s N L \ee^{-t \Delta}-\Tr_s N L \ee^{-t \Delta_1}-\Tr_s N L \ee^{-t \Delta_1}) dt\\
&=&-2(2T_1(M)-T_1(\tilde{M}_1)-T_1(\tilde{M}_2))\ .\end{eqnarray*}
This finishes the proof of the theorem. \hB

\section{Products and coverings}\label{covprod}

Let $G$ be a compact connected Lie group and $\Gamma$ be a finite group.
Let $C(\Gamma)$ denote the algebra of $\C$-valued functions on $\Gamma$.
We need the generalization of higher equivariant analytic torsion
$T^\Gamma(M)\in \tilde{I}(G)\otimes C(\Gamma)$ mentioned in \cite{lott94}, Note 3.
Let $M$ be a closed oriented $G\times \Gamma$-
manifold equipped with a $G\times \Gamma$-invariant Riemannian metric
$g^M$. Set 
$$\chi^\prime(M)(\gamma):=\sum_{p=0}^\infty p (-1)^p \Tr \: H^p(\gamma)\ ,$$
where $H^p(\gamma)$ is the induced action of $\gamma\in\Gamma$ on $H^p(M,\R)$.
Then we define
$T^\Gamma(M)\in \tilde{I}(G)\otimes C(\Gamma)$ to be the element
represented by the function
$$\gamma\mapsto \frac{d}{ds}_{|s=0}\frac{1}{\Gamma(s)}\int_0^\infty (\Tr_s N \gamma \ee^{D_t^2} t^{s-1}-\chi^\prime(M)(\gamma)) dt\ .$$

Let $M$ be a closed oriented $G\times \Gamma$-manifold and
$N$ be a closed oriented $\Gamma$-manifold.
Then we form the closed oriented $G\times \Gamma$-manifold $M\times N$, where
$\Gamma$ acts diagonally. We choose a $G\times\Gamma$-invariant
Riemannian metric $g^M$, a $\Gamma$-invariant Riemannian metric
$g^N$, and we let $g^{M\times N}$ be the product metric.

Define the $\Gamma$-equivariant Euler characteristic
$\chi^\Gamma(N)\in C(\Gamma)$ of a closed $\Gamma$-manifold $N$ by
$$\chi^\Gamma(N)(\gamma):=\sum_{p=0}^\infty (-1)^p \Tr\:H^p(\gamma)\ .$$
\begin{lem}\label{prdf}
If $\chi^\Gamma(M)=0$, then
$$T^\Gamma(M\times N)=T^\Gamma(M)\chi^\Gamma(N)\ .$$
\end{lem}
\proof
We write $D_t(M), D_t(N),D_t(M\times N$ for the operators (\ref{hamil}) on $M,N,M\times N$.
Let $\Delta(N)$ be the Laplace operator on $\Omega(N)$.
On the level of Hilbert space closures we have
$$\clo_{L^2}\Omega(M\times N)=\clo_{L^2}\Omega(M)\otimes\clo_{L_2}\Omega(N)\ .$$
With respect to this splitting we can write
$$D_t(M\times N)^2=D_t(M)^2\times 1 - 1\otimes t\Delta(N)\ .$$
If $\gamma\in\Gamma$, then
\begin{eqnarray*}
&&\Tr_s N\gamma \ee^{D_t(M\times N)^2}\\
&=&\Tr_s (N\otimes 1 + 1\otimes N)(\gamma\otimes\gamma) (\ee^{D_t(M)^2}\otimes\ee^{-t\Delta(N)})\\
&=&\Tr_s N \gamma \ee^{D_t(M)^2} \Tr_s \gamma \ee^{-t\Delta(N)} + \Tr_s \gamma \ee^{D_t(M)^2} \Tr_s N \gamma \ee^{-t\Delta(N)}\ .
\end{eqnarray*}
By the equivariant McKean-Singer formula \cite{berlinegetzlervergne92}, Thm. 6.3, we have $\Tr_s \gamma \ee^{-t\Delta(N)}=\chi^\Gamma(N)(\gamma)$.
Moreover we have
\begin{eqnarray*}
\frac{d}{dt}\Tr_s \gamma \ee^{D_t(M)^2}
&=&\Tr_s \gamma \frac{d}{dt}D_t^2 \ee^{D_t(M)^2}\\
&=&\Tr_s [ \frac{d}{dt}D_t, \gamma D_t\ee^{D_t(M)^2}]\\
&=&0\\
\lim_{t\to \infty }\Tr_s \gamma \ee^{D_t(M)^2}&=&
\lim_{t\to \infty }\Tr_s \gamma \ee^{-t\Delta(M)}\\
&=&\chi^\Gamma(M)(\gamma)\\
&=& 0\ .
\end{eqnarray*}
It follows
$$\Tr_s N\gamma \ee^{D_t(M\times N)^2}= \chi^\Gamma(N)(\gamma) \Tr_s N \gamma \ee^{D_t(M)^2}\ .$$
This implies the assertion of the Lemma. \hB

Let $N$ be a closed oriented $G\times\Gamma$-manifold such that
$\Gamma$ acts freely on $N$. Let $M:=\Gamma\backslash N$.
Then $M$ is a closed oriented $G$-manifold.
We equip $N$ with a $G\times \Gamma$-invariant Riemannian metric
and define  $g^M$ such that the projection $\pi:N\rightarrow M$ becomes
a local isometry.

Let $\int_\Gamma:C(\Gamma)\rightarrow \C$ be the integral over $\Gamma$
with respect to the normalized Haar measure.
We denote the induced map $\tilde{I}(G)\otimes C(\Gamma)\rightarrow \tilde{I}(G)$
by the same symbol.

\begin{lem}\label{covform}
$$T(M)=\int_\Gamma \:T^\Gamma(N)\ .$$
\end{lem}
\proof
Note that $\Pi:=\frac{1}{|\Gamma|}\sum_{\gamma\in \Gamma} \gamma$
acts on $\Omega(N)$ as projection onto the subspace of  $\Gamma$-invariant
forms which can be identified with $\Omega(M)$ using the pull-back $\pi^*$.
Moreover, $D_t(M)$ coincides with the restriction of $D_t(N)$ to the range
of  $\Pi$.
We have
\begin{eqnarray*}
\frac{1}{|\Gamma|}\sum_{\gamma\in \Gamma} \Tr_s N\gamma \ee^{D_t(N)^2}&=&
\Tr_s N \Pi \ee^{D_t(N)^2}\\
&=&\Tr_s N \ee^{D_t(M)^2}\ .
\end{eqnarray*}
This implies the  assertion of the Lemma. \hB

\section{Manifolds with corner singularities}\label{formalities}

In this section we extend the definition of $T(M)$, $T^\Gamma(M)$,
and the results of Section \ref{covprod} to manifolds with
corner singularities.

A compact manifold with a corner singularitiy of codimension one
is just a manifold with boundary. Corner singularities of
codimension two arise if we admit that boundaries have itself
boundaries. In general a corner singularity of codimension $m$ of a $n$-dimensional
manifold is modelled on $(R_+)^m\times R^{n-m}$, where $R_+=[0,\infty)$.

Let $M$ be a compact manifold with corner singularities.
Then the boundary of $M$ can be decomposed into pieces $\partial_1 M\cup\dots\cup\partial_l M$.
We do not require that the pieces $\partial_i M$ are
connected. If $x\in M$ belongs to a corner singularity of codimension $m$,
then $x$ meets exactly $m$ pieces of $\partial M$.

For $i\in \{1,\dots l\}$ we can form
the double $\tilde{M}_i:=M\cup_{\partial_i M}M$ of $M$ along the piece $\partial_i M$.
Then $\tilde{M}_i$ is again a compact manifold with corner singularities.
In particular it has $l-1$ boundary pieces $\partial_j\tilde{M}_i=\partial_j M\cup_{\partial_j M\cap \partial_i M}\partial_j M$, $j\not=i$.

The notion corner singularities and the construction of the double
extends to compact oriented $G$-manifolds in the obvious way.
We define $T(M)$ for compact oriented $G$-manifolds inductively with
respect to the number $l(M)$ of boundary pieces.

If $l(M)=0$, then $T(M)$ is  already defined.
Assume now that $l(M)$ is defined for all $M$ with $l(M)<l$.
If $M$ is now a compact oriented $G$-manifold with $l(M)=l$.
Then we set
$$T(M):=\frac{1}{2} T(\tilde{M_1})\ .$$
If $l>1$, then we have to check that this definition is independent of
the numbering of boundary components.
It suffices to show that $T(\tilde{M_1})=T(\tilde{M_2})$.
Note that  
$\tilde{\tilde{M}_1}_2$ and $\tilde{\tilde{M}_2}_1$   are $G$-diffeomorphic.
Using the induction hypothesis
$$2T(\tilde{M_1})=T(\tilde{\tilde{M}_1}_2)=T(\tilde{\tilde{M}_2}_1)=2T(\tilde{M_2})\ .$$ 
Thus $T(M)$ is well defined.

The doubling trick was introduced by
\cite{lottrothenberg91},  Ch. IX. 
Instead of the formal definition above
one could also employ absolute and relative boundary
boundary conditions in order to define higher equivariant analytic
torsion $T(M,abs)$, $T(M,rel)$ for $G$-manifolds with boundary.
If the Riemannian metric is choosen to be product near
the boundary, then $T(M)=\frac{1}{2} T(M,abs)+T(M,rel)$.

The sum formula \ref{sumff} has now the nice reformulation
\begin{equation}\label{sumfff}
T(M)=T(M_1)+T(M_2) \ .
\end{equation}
It has the following generalization:
\begin{kor}\label{summ2}
Let $M_i$, $i=1,2$, be compact oriented $G$-manifolds with corner singularities.
If we are given a $G$-diffeomorphism $\partial_1 M_1\cong \partial_1 M_2$,
then we form the manifold with corner singularities 
$M:=M_1\cup_{\partial_1 M_i}M_2$, and we have
$$T(M)=T(M_1)+T(M_2)\ .$$
\end{kor}
\proof
We employ induction by the number of boundary pieces.
The assertion is true if $M$ is closed.
Assume that the corollary holds true for all $M$ with $l(M)<l$.
Let $M=M_1\cup_{\partial_1 M_i}M_2$ now be a manifold with $l(M)=l$ and $l\ge 1$. Then we can assume that $l(M_1)\ge 2$.
Let $\partial_1 M$ be the piece corresponding to $\partial_2 M_1$.
We distinguish the cases (a):  $\partial_2 M_1\cap \partial 1M_1=\emptyset$
and (b): $\partial_2 M_1\cap \partial_1M_1\not=\emptyset$.
In case (a) let $\partial_1 M$ be the piece corresponding to $\partial_2 M_1$.
Then using the induction hypothesis
$$T(M)=\frac{1}{2}T(\tilde{M}_1)=\frac{1}{2}T(\tilde{M_1}_2)+T(M_2)=T(M_1)+T(M_2)\ .$$
In case (b) there is a boundary piece $\partial_2 M_2$ meeting $\partial_1 M_2$.
Then $M$ has a boundary piece $\partial_1 M:=\partial_2 M_1\cup_{\partial_1 M_i\cap \partial_2 M_i}\partial_2 M_2$.
Again using the induction hypothesis we have
$$T(M)=\frac{1}{2}T(\tilde{M}_1)=\frac{1}{2}T(\tilde{M_1}_2) +  \frac{1}{2}T(\tilde{M_2}_2)=T(M_1)+T(M_2)\ .$$
This proves the corollary. \hB

Let $\Gamma$ be an additional finite group.
For a $G\times\Gamma$-manifold with corner singularities
we require that that the pieces $\partial_i M$ are
compact $G\times \Gamma$-manifolds with corner singularities as well.

A Riemannian metric on a manifold with corner singularities is compatible
if it is a product metric $g^{(R_+)^m} + g^{\R^{n-m}}$ at a corner of codimension $m$.
Then we can form the doubles $\tilde{M}_i$ metrically.

Let $M$ be a compact oriented $G\times \Gamma$-manifold
with corner singularities equipped with a compatible $G\times\Gamma$-
invariant Riemannian metric.
Then we define
$T^\Gamma(M)$ for $G\times\Gamma$-manifolds with corner singularities
using the same formal procedure as for trivial $\Gamma$.
We can generalize Lemma \ref{covform} to this case.
Let $N$ be a compact oriented $G\times\Gamma$-manifold with corner singularities
such that $\Gamma$ acts freely and form $M:=\Gamma\backslash N$.
\begin{kor}\label{covformex}
$$T(M)=\int_\Gamma T^\Gamma(N)\ .$$
\end{kor}
\proof
We argue by induction with respect to the number of boundary pieces.
If $l(N)=0$, then this is just Lemma \ref{covform}.
Assume now that the corollary holds true for all $N$ with $l(N)<l$.
Let now $N$ be a compact oriented $G\times\Gamma$-manifold with corner singularities
such that $\Gamma$ acts freely and $l(N)=l\ge 1$.
Then consider the covering $\tilde{N}_1\rightarrow \tilde{M}_1$.
Applying the induction hypothesis we obtain
$$T(M)=\frac{1}{2} T(\tilde{M}_1)=\frac{1}{2}\int_\Gamma T^\Gamma(\tilde{N}_1)=\int_\Gamma T^\Gamma(N)\ .$$
This proves the corollary. \hB

Let $M$ be a closed oriented $G\times \Gamma$-manifold
and $N$ be a compact oriented $\Gamma$-manifold with corner singularities.
Then we form the compact oriented $G\times \Gamma$-manifold $M\times N$ with
corner singularities, where $\Gamma$ acts diagonally.
We choose a $G\times\Gamma$-invariant
Riemannian metric $g^M$, a $\Gamma$-invariant compatible Riemannian metric
$g^N$, and we let $g^{M\times N}$ be the product metric which is again invariant and compatible.

We define the $\Gamma$-equivariant Euler characteristic
$\chi^\Gamma(N)\in C(\Gamma)$ of a $\Gamma$-manifold $N$
with corner singularities with $l(N)\ge 1$
inductively with respect to the number of boundary pieces by
$$\chi^\Gamma(N):=\frac{1}{2}\chi^\Gamma(\tilde{N}_1)\ .$$
We leave it to the interested reader to express $\chi^\Gamma(N)$
in terms of equivariant Euler characteristics of the components of the
filtration of $N$. The main feature of this definition is that
the  equivariant Euler characteristic is additive under glueing
along boundary pieces.

We have the following generalization of Lemma \ref{prdf}.
\begin{kor}\label{prdfex}
If $\chi^\Gamma(M)=0$, then
$$T^\Gamma(M\times N)=T^\Gamma(M)\chi^\Gamma(N)\ .$$
\end{kor}
\proof
We argue by induction over the number of boundary  pieces $l(N)$.
If $l(N)=0$, then this is just Lemma \ref{prdf}.
Assume that the corollary holds true if $l(N)<l$.
Let now $N$ be such that $l(N)=l\ge 1$.
Let $\partial_1(M\times N):=M\times \partial_1 N$.
Then using the induction hypothesis and the additivity of $\chi^\Gamma$
we obtain
$$T^\Gamma(M\times N)=\frac{1}{2}T^\Gamma((\widetilde{M\times N})_1)=\frac{1}{2}T^\Gamma(M)\chi^\Gamma(\tilde{N}_1)=
T^\Gamma(M)\chi^\Gamma(N)\ .$$
This proves the corollary. \hB

\section{Restriction to subgroups}\label{subgr}

Let $G,H$ be a  connected compact Lie groups with Lie algebras $g,h$.
If $f:H\rightarrow G$ is a homomorphism, then $f_*:h\rightarrow g$
induces a map $\tilde{f}:I(G)\rightarrow I(H)$.
If $H\subset G$ is a closed subgroup and $i$ denotes the inclusion,
then we set $\tilde{i}=:\res^G_H$.
If $g\in G$, then we put $H^g:=gHg^{-1}$. Let $\alpha_g:H\rightarrow H^g$
be given by $\alpha_g(h):=ghg^{-1}$.

Let $M$ is a closed oriented $G$-manifold with corner singularities.
If $f:H\rightarrow G$ is a homomorphism, then we denote
by $f^* M$ the $H$-manifold $M$ with action induced by $f$. 
If $H\subset G$ is a closed subgroup, then we put
$\res^G_H M:=i^* M$.
The following Lemma is  an immediate consequence of the definition
of  $T(M)$.
\begin{lem}\label{ress}
{\bf (1)$\:\:$:} If $f:H\rightarrow G$ is a homomorphism, then $\tilde{f}T(M)=T(f^*M)$. In particular, if $H\subset G$ is closed,
then  $\res^G_H T(M)=T(\res^G_H M)$.\\
{\bf (2)$\:\:$:}If $H\subset G$ is closed, then for all $g\in G$ we have $\tilde{\alpha}_g\res^G_H T(M) = \res^G_{H^g} T(M)$ for all $g\in G$.
\end{lem}
The association $H\subset G\mapsto \tilde{I}(H)=:\tilde{I}_G(H)$ assembles to give a
contravariant functor $\tilde{I}_G:\Or(G)\rightarrow \C-vect$.
If $f:H\rightarrow G$ is a homomorphism, then it induces
a natural functor $f_*:\Or(H)\rightarrow \Or(G)$ sending
$H/K$ to $G/f(K)$. For $K\subset H$ let $f_K:K\rightarrow f(K)$ be the
restriction of $f$ to $K$. The collection $\{\tilde{f}_K\}$, $K\in H$,
provides a natural transformation $\tilde{f}:\tilde{I}_G\circ f_*\rightarrow \tilde{I}_H$.
Let $f^*:\lim_{\Or(G)} \tilde{I}_G \rightarrow \lim_{\Or(H)}\tilde{I}_H$
denote the induced map.

Lemma \ref{ress} says that $G/H\mapsto T(\res^G_H M)$ is a section
of $\tilde{I}_G$. Since $\Or(G)$ has a final object $G/G$,
we have an isomorphism
\begin{equation}\label{ident}\lim_{\Or(G)} \tilde{I}_G\cong \tilde{I}(G)\end{equation}
given by restriction to the final object.

By $S(G)$ we denote the full  subcategory of $\Or(G)$
of those objects $G/H$ with $H\cong S^1$.
We denote the space of sections of $\tilde{I}_{G|S(G)}$ by $V(G)$, i.e.
$$V(G):=\lim_{S(G)} \tilde{I}_G\ .$$
There is a natural restriction map
$$R_G:\tilde{I}(G)\cong\lim_{\Or(G)} \tilde{I} \rightarrow V(G)\ .$$
\begin{lem} \label{ress1}
$R_G$ is injective.
\end{lem}
\proof
Let $T\subset G$ be a maximal torus and denote by $j$ its inclusion.
There is a functor $j_{*|S(T)}:S(T)\rightarrow S(G)$. Let
$J^*:\lim_{S(G)} \tilde{I}_G \rightarrow \lim_{S(T)} \tilde{I}_T$
be induced by the natural tranformation
$\tilde{j_{*|S(T)}}:\tilde{I}_{G|S(G)}\circ j_{*|S(T)}\rightarrow \tilde{I}_{T|S(T)}$.
Then $R_T\circ j^*=J^*\circ R_G$.
In order to prove that $R_G$ is injective it is therefore
sufficient to show that $j^*$ and $R_T$ are injective.

Now $j^*$ is injective  since it coincides with
$\res^G_T :\tilde{I}(G)\rightarrow \tilde{I}(T)$ under the identification
(\ref{ident}), and the latter map well known to be injective.
Let $t$ be the Lie algebra of $T$. The kernel of $\exp:t\rightarrow T$
defines a $\Z$ structure on $t$. The set of subspaces
$h\subset t$ corresponding
to objects $T/H\in S(T)$ with $H\cong S^1$ is just
the set of integral points of the projective space $P(t\otimes \C)$.
Injectivity of $R_T$ follows easily from the fact
that the set of integral points of $P(t\otimes \C)$ is Zariski dense. \hB

\begin{kor}\label{fgd}
$T(M)$ is uniquely determined by the values of $T(\res^G_H M)$
for all $H\subset G$ with $H\cong S^1$.
\end{kor}

\section{The map $T_G$}\label{s1}

We need the following technical result.
\begin{lem}\label{decomp}
Let $M$ be a closed manifold. Then there exists a Riemannian
metric $g^M$ and a decomposition $M=\cup_i B_i$ of $M$ into
manifolds with corner singularities such that the $B_i$ are contractible 
and the restriction of
$g^M$ to $B_i$ is compatible for all $i$.
\end{lem}
\proof
We choose a smooth triangulation of $M$. Then there is another
smooth triangulation $\cT$ which is dual to the first one.
We choose small closed tubular neighbourhoods $U_\sigma$ of the simplices $\sigma$ of $\cT$. We now proceed inductively. Assume that in the steps $0,\dots,l-1$
we have already  defined $B_i$, $i=1,\dots r$.
In the $l$'th step we let $B_{r+1},\dots$ be the
intersections of $U_\sigma\cap (M\setminus \inter(\cup_{i=1}^r B_i))$,
where $\sigma$ runs over all simplices of $\cT$ of dimension $j$.
By choosing the tubular neighbourhoods appropriately, this 
construction gives manifolds $B_i$ with corner singularities.
Now one can construct an appropriate Riemannian metric.
\hB

Recall that if  $M$ is a manifold with corner singularities and $M$ has
at least one boundary piece, then we define inductively
$\chi(M):=\frac{1}{2}\chi(\tilde{M}_1)$. In particular,
if $M=\cup_i B_i$ is a decomposition as in Lemma \ref{decomp}, then
\begin{equation}\
\label{p1oi} \chi(M)=\sum_i \chi(B_i)\ .\end{equation}

Recall the definition (\ref{ts1}) of $T_{S^1}:U(S^1)\rightarrow \tilde{I}(S^1)$.
\begin{prop}
{\bf (1)$\:\:$:} Let $M$ be a closed oriented $S^1$-manifold.
Then $T(M)=T_{S^1}\chi_{S^1}(M)$.\\
{\bf (2)$\:\:$:} If in addition $M$ is even-dimensional, then $\chi_{S^1}(M)=a[S^1/S^1]$ for some $a\in \Z$.
\end{prop}
\proof
A compact $S^1$-manifold has a finite number of orbit types
$H_1,\dots,H_l$. We employ induction by the number of orbit types $l(M)$.
We first assume that $l(M)=1$. If $H_1=S^1$, then
$T(M)=0$ and $\chi_{S^1}(M)=\chi(M)[S^1/S^1]$.
Thus $T_{S^1}\chi_{S^1}(M)=\chi(M)T_{S^1}[S^1/S^1]=0$,
too.

We now consider the case that  $H_1\not=S^1$.
Then by \cite{bredon72}, II.5.2., we have a smooth locally trivial fibre bundle $M\rightarrow M/S^1$
with fibre $S^1/H_1$.
Let $M/S^1=\cup_i B_i$ be a decomposition of $M/S^1$
into manifolds  with corner  singularities given by Lemma \ref{decomp}.
Then $M_{|B_i}\cong S^1/H_1\times B_i$.
Using Corollaries \ref{summ2} and \ref{prdfex},
 (\ref{ts1}),  (\ref{p1oi}), and 
$\chi_{S^1}(M)=\chi(M/S^1)[S^1/H_1]$  
we obtain
\begin{eqnarray*}
T(M)&=&\sum_i T(M_{|B_i})\\
&=&\sum_i T(S^1/H_1\times B_i)\\
&=&T(S^1/H_1) \sum_i \chi(B_i)\\
&=&T_{S^1}[S^1/H_1]\chi(M/S^1)\\
&=&T_{S^1} \chi_{S^1}(M)\ .
\end{eqnarray*}
This proves assertion {\bf (1)} for $l(M)=1$.
If $M$ is even-dimensional closed,
then $M/S^1$ is odd-dimensional, and $\chi(M/S^1)=0$ by Poincar\'e
duality. Assertion {\bf (2)} follows.

Now assume that the proposition  holds true for all
$M$ with $l(M)<l$.
Let $M$ be a closed oriented $S^1$-manifold with $l(M)=l$.
Without loss of generality we can assume that $H:=H_1\not=S^1$.
By \cite{bredon72}, VI 2.5., the fixed point set $M_H$ of $H$ is a smooth
submanifold of $M$ with normal bundle $NM_H$,
which we identify with an equivariant
tubular neighbourhood of $M_H$
using the exponential map provided by a $S^1$-invariant Riemannian metric $g^M$.

Assume that $M$ is odd-dimensional.
By Corollary \ref{summ2} we have
$T(M)=T(M\setminus NM_H) + T(\bar{NM_H})$.
Let $N$ be the double of $M\setminus NM_H$.
Then $l(N)\le l-1$, and we can apply the induction hypothesis
in order to obtain
$T(M\setminus NM_H)=\frac{1}{2}T(N)=\frac{1}{2}T_{S^1}\chi_{S^1}T(N)$.
Note that $$\chi_{S^1}(N)=2\chi_{S^1}(M\setminus NM_H)-\chi_{S^1}(\partial \bar{NM_H})\ .$$
Note that
$\partial \bar{NM_H}$ is even-dimensional, closed and orientable. Since
$l(\partial \bar{NM_H})<l$ we have by our induction hypothesis
$\chi_{S^1}(\partial \bar{NM_H})=a[S^1/S^1]$ for some $a\in\Z$. This implies
$T_{S^1}\chi_{S^1}(\partial \bar{NM_H})=0$ and
\begin{equation}\label{teil1}
T(M\setminus NM_H)=T_{S^1}\chi_{S^1}(M\setminus NM_H)\ .
\end{equation}

We now compute $T(\bar{NM_H})$.
Since $l(M_H)=1$ we have a smooth locally trivial fibre bundle
$M_H\rightarrow M_H/S^1$ with fibre $S^1/H$.
Let $M_H/S^1=\cup_i B_i$ be a decomposition of $M_H/S^1$ into
manifolds with corner singularities given by Lemma \ref{decomp}.
Then
$M_{H,i}:=(M_H)_{|B_i}\cong S^1/H\times B_i$.
Since $H$ acts orientation preserving, the bundle $NM_H$ admits an $H$-invariant
complex structure. 
The restriction $(NM_H)_{|M_{H,i}}$ can be written  as
$S^1\times V_i/H$, where $V_i\rightarrow B_i$ is a complex vector bundle
on which $H$ acts fibrewise linear.

Since a complex linear action of a cyclic group $H$
can always be  extended to the connected group $S^1$,
we obtain $\chi^H(V_i)(\gamma)=\chi(V_i)$ for all $\gamma\in H$.
Moreover we have $\chi^H(S^1)=0$. Thus we can apply
Corollaries \ref{prdfex} and \ref{covformex} in order to obtain
$$T(S^1\times V_i/H)=\int_\Gamma T^H(S^1)\chi^H(V_i)=\int_\Gamma T^H(S^1)\chi(V_i)=T(S^1/H)\chi(V_i)\ .$$
Since $\bar{NM_H}$ and $M_H$
are $S^1$-homotopy equivalent, we have
$\chi_{S^1}(\bar{NM_H})=\chi_{S^1}(M_H)$.
Moreover, $\sum_i\chi(V_i)=\sum_i\chi(B_i)=\chi(M_H/S^1)$
and $\chi_{S^1}(M_H)=\chi(M_H/S^1) [S^1/H]$.
Thus we obtain by Corollary \ref{summ2}
\begin{eqnarray}
T(\bar{NM_H})&=&\sum_i T(S^1\times V_i/H)\nonumber\\
&=&\sum_i T(S^1/H)\chi(V_i)\nonumber\\
&=&T_{S^1}[S^1/H] \chi(M_H/S^1)\nonumber\\
&=&T_{S^1}\chi_{S^1}(\bar{NM_H})\label{rt2}.
\end{eqnarray}
We have
$$\chi_{S^1}(M)=\chi_{S^1}(M\setminus NM_H)+\chi_{S^1}(\bar{NM_H})-\chi_{S^1}(\partial\bar{NM_H})\ .$$
Since $T_{S^1} \chi_{S^1}(\partial\bar{NM_H})=0$, combining (\ref{teil1}) and (\ref{rt2})
we obtain the desired formula $T(M)=T_{S^1}\chi_{S^1}(M)$ for $M$ odd-dimensional.

Assume now that $M$ is even-dimensional and that $l(M)=l$.
Then $T(M)=0$,  and {\bf (1)} follows from {\bf (2)}.
We now show {\bf (2)}. 
We have
$$\chi_{S^1}(M)=\chi_{S^1}(M\setminus NM_H)+\chi_{S^1}(M_H)\ .$$
We can apply the induction hypothesis to $M_H$
and the double of $M\setminus NM_H$. It follows that
$\chi_{S^1}(M\setminus NM_H)=\frac{1}{2}\chi_{S^1}(\partial \bar{NM_H})+
a[S^1/S^1]$.
The restriction $\partial \bar{NM_H}_{|M_{H,i}}$
is isomorphic to  
$S^1\times \partial  S V_i/H$, where $S V_i$ denotes the sphere bundle of $V_i$.

Let $U$ be the unit sphere in a fibre of $NM_H$.  
Using that $M_H/S^1$ is closed, orientable, and
odd-dimensional,  we obtain
\begin{eqnarray*}
\chi_{S^1}(\partial \bar{NM_H})&=&
\sum_i \chi_{S^1}(S^1\times U/H)\chi(B_i)\\
&=&\chi_{S^1}(S^1\times U/H) \chi(M_H/S^1)\\
&=&0\ .
\end{eqnarray*}
This finishes the proof of {\bf (2)}. \hB

We now construct $T_G$. The collection $T_H$, $H\in S(G)$,
forms a natural transformation from the functor $H\mapsto U(H)$
to $H\mapsto \tilde{I}(H)$. Thus we obtain a homomorphism
$$\tilde{T}:\lim_{S(G)} U\rightarrow V(G)\ .$$

Let $\tilde{\res} : U(G)\rightarrow  \lim_{S(G)}$ be given by
the collection $\res^G_H$, $H\in S(G)$.
If $M$ is a compact $G$-manifold, then we let $\tilde{\chi}(M)\in \lim_{S(G)}U$
be given by the section $S(G)\ni H\mapsto \chi_H(M)\in U(H)$.
Then $\tilde{\res}  \chi_G(M)=\tilde{\chi}(M)$.
\begin{lem}\label{tge}
There is a unique homomorphism $T_G:U(G)\rightarrow \tilde{I}(G)$
such that $R_G\circ T_G=\tilde{T}\circ  \tilde{\res}$.
\end{lem}
\proof
For $G/K\in \Or(G)$ we shall have
\begin{eqnarray*}
R_G\circ T_G[G/K]&=&\tilde{T}\circ \tilde{\res} \chi_G(G/K)\\
&=&\tilde{T}\circ\tilde{\chi}(G/K)\\
&=&\{S(G)\ni H\mapsto T_H \circ\chi_H \circ\res^G_H(G/K)\}\\
&=&\{S(G)\ni H\mapsto T(\res^G_H\:G/K)\}\\
&=&\{S(G)\ni H\mapsto \res^G_H T(G/K)\}\\
&=&R_G T(G/K) \ .
\end{eqnarray*}
Hence by injectivity of $R_G$ (Lemma \ref{ress1})
we are forced to define $T_G[G/K]:=T(G/K)$. \hB

We now finish the proof of Theorem \ref{main}.
Let $M$ be a closed oriented $G$-manifold. Then we have
\begin{eqnarray*}
R_G\circ T_G\chi_G(M)&=&\tilde{T}\circ \tilde{\res} \chi_G(M)\\
&=&\tilde{T}\circ \tilde{\chi}(M)\\
&=&\{S(G)\ni H\mapsto T_H \circ\chi_H \circ\res^G_H(M)\}\\
&=&\{S(G)\ni H\mapsto T(\res^G_H\:M)\}\\
&=&\{S(G)\ni H\mapsto \res^G_H T(M)\}\\
&=&R_G T(M) \ .
\end{eqnarray*}
We conclude that $T_G\chi_G(M)=T(M)$ by Lemma \ref{ress1}.\hB

\bibliographystyle{plain}

\end{document}